\documentclass[reqno]{amsart}
\usepackage[square]{natbib}
\newtheorem{thm}{Theorem}
\newtheorem{lem}[thm]{Lemma}
\newtheorem{cor}[thm]{Corollary}

\theoremstyle{definition}

\newtheorem{rem}[thm]{Remark}

\providecommand{\abs}[1]{\lvert#1\rvert}
\providecommand{\Abs}[1]{\Bigl\lvert#1\Bigr\rvert}

\begin{document}
\title[Randomly reinforced urns]{A central limit theorem and its applications to\\multicolor randomly reinforced urns}
\author{Patrizia Berti}
\address{Patrizia Berti, Dipartimento di Matematica Pura ed Applicata ''G. Vitali'', Universita' di Modena e Reggio-Emilia, via Campi 213/B, 41100 Modena, Italy}
\email{patrizia.berti@unimore.it}
\author{Irene Crimaldi}
\address{Irene Crimaldi, Dipartimento di Matematica, Universita' di
Bologna, Piazza di Porta San Donato 5, 40126 Bologna, Italy}
\email{crimaldi@dm.unibo.it}
\author{Luca Pratelli}
\address{Luca Pratelli, Accademia Navale, viale Italia 72, 57100 Livorno,
Italy} \email{pratel@mail.dm.unipi.it}
\author{Pietro Rigo}
\address{Pietro Rigo (corresponding author), Dipartimento di Economia Politica e Metodi Quantitativi, Universita' di Pavia, via S. Felice 5, 27100 Pavia, Italy}
\email{prigo@eco.unipv.it} 

\date{April 24, 2009. First version: April 6, 2009. 
http://arxiv.org/abs/0904.0932v1}

\begin{abstract} Let $(X_n)$ be a sequence of integrable real random variables, adapted to a filtration $(\mathcal{G}_n)$.
Define
\begin{equation*}
C_n=\sqrt{n}\,\bigl\{\frac{1}{n}
\sum_{k=1}^nX_k-E(X_{n+1}\mid\mathcal{G}_n)\bigr\}\quad\text{and}\quad
D_n=\sqrt{n}\,\bigl\{E(X_{n+1}\mid\mathcal{G}_n)-Z\bigr\}
\end{equation*}
where $Z$ is the a.s. limit of $E(X_{n+1}\mid\mathcal{G}_n)$
(assumed to exist). Conditions for
$(C_n,D_n)\longrightarrow\mathcal{N}(0,U)\times\mathcal{N}(0,V)$
stably are given, where $U,\,V$ are certain random variables. In
particular, under such conditions, one obtains
\begin{equation*}
\sqrt{n}\,\bigl\{\frac{1}{n}
\sum_{k=1}^nX_k-Z\bigr\}=C_n+D_n\longrightarrow\mathcal{N}(0,U+V)\quad\text{stably}.
\end{equation*}This CLT has natural applications to Bayesian
statistics and urn problems. The latter are investigated, by paying
special attention to multicolor randomly reinforced generalized
Polya urns.
\end{abstract}

\maketitle

\section{Introduction and motivations}\label{intro}

As regards asymptotics in urn models, there is not a unique
reference framework. Rather, there are many (ingenious) disjoint
ideas, one for each class of problems. Well known examples are
martingale methods, exchangeability, branching processes, stochastic
approximation, dynamical systems and so on; see \cite{P}.

Those limit theorems which unify various urn problems, thus, look of
some interest.

In this paper, we focus on the CLT. While thought for urn problems,
our CLT is stated for an arbitrary sequence $(X_n)$ of real random
variables. Accordingly, it potentially applies to every urn
situation, but it has generally a broader scope. Suppose
$E\abs{X_n}<\infty$ and define
$Z_n=E\bigl(X_{n+1}\mid\mathcal{G}_n\bigr)$ where $(\mathcal{G}_n)$
is some filtration which makes $(X_n)$ adapted. Under various
assumptions, one obtains $Z_n\overset{a.s.,L_1}\longrightarrow Z$
for some random variable $Z$. Define further
$\overline{X}_n=\frac{1}{n} \sum_{k=1}^nX_k$ and
\begin{gather*}
C_n=\sqrt{n}\,\bigl(\overline{X}_n-Z_n),\quad
D_n=\sqrt{n}\,\bigl(Z_n-Z),\quad
W_n=\sqrt{n}\,\bigl(\overline{X}_n-Z).
\end{gather*}

The limit distribution of $C_n$, $D_n$ or $W_n$ is a main goal in
various fields, including Bayesian statistics, discrete time
filtering, gambling and urn problems. See \cite{AMS}, \cite{BCL},
\cite{BPR}, \cite{BCPR}, \cite{CLP}, \cite{C}, \cite{GR} and
references therein. In fact, suppose the next observation $X_{n+1}$
is to be predicted basing on the available information
$\mathcal{G}_n$. If the predictor $Z_n$ cannot be evaluated in
closed form, one needs some estimate $\widehat{Z}_n$ and $C_n$
reduces to the scaled error when $\widehat{Z}_n=\overline{X}_n$. And
$\overline{X}_n$ is a sound estimate of $Z_n$ under some
distributional assumptions on $(X_n)$, for instance when $(X_n)$ is
exchangeable, as it is usual in Bayesian statistics. Similarly,
$D_n$ and $W_n$ are of interest provided $Z$ is regarded as a random
parameter. In this case, $Z_n$ is the Bayesian estimate (of $Z$)
under quadratic loss and $\overline{X}_n$ can be often viewed as the
the maximum likelihood estimate. Note also that, in the trivial case
where $(X_n)$ is i.i.d. and $\mathcal{G}_n=\sigma(X_1,\ldots,X_n)$,
one obtains $C_n=W_n=\sqrt{n}\,\bigl(\overline{X}_n-EX_1)$ and
$D_n=0$. As to urn problems, $X_n$ could be the indicator of
$\{$black ball at time $n\}$ in a multicolor urn. Then, $Z_n$
becomes the proportion of black balls in the urn at time $n$ and
$\overline{X}_n$ the observed frequency of black balls at time $n$.

Our main result (Theorem \ref{main}) provides conditions for
\begin{equation}\label{piot} (C_n,D_n)\longrightarrow
\mathcal{N}(0,U)\times \mathcal{N}(0,V)\quad\text{stably}
\end{equation}
where $U,\,V$ are certain random variables and $\mathcal{N}(0,L)$ is
the Gaussian kernel with mean $0$ and variance $L$. A nice
consequence is that
\begin{equation*} W_n=C_n+D_n\longrightarrow
\mathcal{N}(0,U+V)\quad\text{stably}.
\end{equation*}
Stable convergence, in the sense of Aldous and Renyi, is a strong form of
convergence in distribution; the definition is recalled in Section
\ref{prel}.

To check the conditions for \eqref{piot}, it is fundamental to know
something about the convergence rate of
\begin{gather*}
Z_{n+1}-Z_n=E\bigl(X_{n+2}\mid\mathcal{G}_{n+1}\bigr)-E\bigl(X_{n+1}\mid\mathcal{G}_n\bigr),
\\E\bigl(Z_{n+1}-Z_n\mid\mathcal{G}_n\bigr)=E\bigl(X_{n+2}-X_{n+1}\mid\mathcal{G}_n\bigr).
\end{gather*}
If $(X_n)$ is conditionally identically distributed with respect to
$(\mathcal{G}_n)$, in the sense of \cite{BPR}, then $(Z_n)$ is a
$(\mathcal{G}_n)$-martingale and thus only $Z_{n+1}-Z_n$ plays a
role. This happens in particular if $(X_n)$ is exchangeable and
$\mathcal{G}_n=\sigma(X_1,\ldots,X_n)$.

To illustrate how the CLT works, three applications are given:
$r$-step predictions, Poisson-Dirichlet sequences, and {\em randomly
reinforced generalized Polya urns}. We next describe the latter, the
main of such applications, and we refer to Subsections \ref{kstep}
and \ref{poidir} for the remaining two.

An urn contains black and red balls. At each time $n\geq 1$, a ball
is drawn and then replaced together with a random number of balls of
the same color. Say that $B_n$ black balls or $R_n$ red balls are
added to the urn according to whether $X_n=1$ or $X_n=0$, where
$X_n$ is the indicator of $\{$black ball at time $n\}$. Suppose
\begin{gather*}
B_n\geq 0,\quad R_n\geq 0,\quad EB_n=ER_n\quad\text{for all }n,
\\\sup_nE\bigl\{(B_n+R_n)^u\bigr\}<\infty\quad\text{for some }u>2,
\\m:=\lim_nEB_n>0,\quad q:=\lim_nEB_n^2,\quad \quad s:=\lim_nER_n^2.
\end{gather*}
Letting $\mathcal{G}_n=\sigma(X_1,B_1,R_1,\ldots,X_n,B_n,R_n)$,
suppose also that $(B_{n+1},R_{n+1})$ is independent of
$\mathcal{G}_n\vee\sigma(X_{n+1})$. Then, as shown in Corollary
\ref{poi}, the conditions for \eqref{piot} are satisfied with
\begin{equation*}
U=Z(1-Z)\,\bigl(\frac{(1-Z)q+Zs}{m^2}-1\bigr) \quad\text{and}\quad
V=Z(1-Z)\,\frac{(1-Z)q+Zs}{m^2}.
\end{equation*}

Corollary \ref{poi} improves the existing result on this type of
urns, obtained in \cite{AMS}, under two respects. First, Corollary
\ref{poi} implies convergence of the pairs $(C_n,D_n)$ and not only
of $D_n$. Hence, one also gets
$W_n\longrightarrow\mathcal{N}(0,U+V)$ stably. Second, unlike
\cite{AMS}, neither the sequence $((B_n,R_n))$ is identically
distributed nor the random variables $B_n+R_n$ have compact support.

By just the same argument used for two color urns, multicolor
versions of Corollary \ref{poi} are easily manufactured. To our
knowledge, results of this type were not available so far. Briefly,
for a $d$-color urn, let $X_{n,j}$ be the indicator of $\{$ball of
color $j$ at time $n\}$ where $n\geq 1$ and $1\leq j\leq d$. Suppose
$A_{n,j}$ balls of color $j$ are added in case $X_{n,j}=1$. The
random variables $A_{n,j}$ are requested exactly the same conditions
asked above to $B_n$ and $R_n$. Then,
\begin{gather*}
\bigl({\bf C_n},\,{\bf D_n}\bigr)\longrightarrow\mathcal{N}_d(0,{\bf
U})\times\mathcal{N}_d(0,{\bf V})\quad\text{stably,}
\end{gather*}
where ${\bf C_n}$ and ${\bf D_n}$ are the vectorial versions of
$C_n$ and $D_n$ while ${\bf U},\,{\bf V}$ are certain random
covariance matrices; see Corollary \ref{prio98poi}.

A last note is the following. In the previous urn, the $n$-th
reinforce matrix is
\begin{equation*}
{\bf A_n}=\text{diag}\bigl(A_{n,1},\ldots,A_{n,d}\bigr).
\end{equation*}
Since $EA_{n,1}=\ldots=EA_{n,d}$, the leading eigenvalue of the mean
matrix $E{\bf A_n}$ has multiplicity greater than 1. Even if
significant for applications, this particular case (the leading
eigenvalue of $E{\bf A_n}$ is not simple) is typically neglected;
see \cite{BH}, \cite{J04}, \cite{J05}, and page 20 of \cite{P}. Our
result, and indeed the result in \cite{AMS}, contribute to fill this
gap.

\section{Stable convergence}\label{prel}

Stable convergence has been introduced by Renyi in \cite{REN} and
subsequently investigated by various authors. In a sense, it is
intermediate between convergence in distribution and convergence in
probability. We recall here basic definitions. For more information,
we refer to \cite{AE}, \cite{CLP}, \cite{HH} and references therein.

Let $(\Omega,\mathcal{A},P)$ be a probability space and $S$ a metric
space. A {\em kernel} on $S$ (or a {\em random probability measure}
on $S$) is a measurable collection $N=\{N(\omega):\omega\in\Omega\}$
of probability measures on the Borel $\sigma$-field on $S$.
Measurability means that
\begin{equation*}
N(\omega)(f)=\int f(x)\,N(\omega)(dx)
\end{equation*}
is $\mathcal{A}$-measurable, as a function of $\omega\in\Omega$, for
each bounded Borel map $f:S\rightarrow\mathbb{R}$.

Let $(Y_n)$ be a sequence of $S$-valued random variables and $N$ a
kernel on $S$. Both $(Y_n)$ and $N$ are defined on
$(\Omega,\mathcal{A},P)$. Say that $Y_n$ converges {\it stably} to
$N$ in case
\begin{gather*}
P\bigl(Y_n\in\cdot\mid H\bigr)\rightarrow E\bigl(N(\cdot)\mid
H\bigr)\quad\text{weakly}
\\\text{for all }H\in\mathcal{A}\text{ such that }P(H)>0.
\end{gather*}
Clearly, if $Y_n\rightarrow N$ stably, then $Y_n$ converges in
distribution to the probability law $E\bigl(N(\cdot)\bigr)$ (just
let $H=\Omega$). Moreover, when $S$ is separable, it is not hard to
see that $Y_n\overset{P}\rightarrow Y$ if and only if $Y_n$
converges stably to the kernel $N=\delta_Y$.

We next mention a strong form of stable convergence, introduced in
\cite{CLP}, to be used later on. Let
$\mathcal{F}_n\subset\mathcal{A}$ be a sub-$\sigma$-field, $n\geq
1$. Say that $Y_n$ converges to $N$ {\em stably in strong sense},
with respect to the sequence $(\mathcal{F}_n)$, in case
\begin{equation*}
E\bigl(f(Y_n)\mid\mathcal{F}_n\bigr)\overset{P}\longrightarrow
N(f)\quad\text{for each }f\in C_b(S)
\end{equation*}
where $C_b(S)$ denotes the set of real bounded continuous functions
on $S$.

Finally, we state a simple but useful fact as a lemma.

\begin{lem}\label{plm}
Suppose that $S$ is a separable metric space and

$C_n$ and $D_n$ are $S$-valued random variables on
$(\Omega,\mathcal{A},P)$, $n\geq 1$;

$M$ and $N$ are kernels on $S$ defined on $(\Omega,\mathcal{A},P)$;

$(\mathcal{G}_n:n\geq 1)$ is an (increasing) filtration satisfying
\begin{equation*}
\sigma(C_n)\subset\mathcal{G}_n\quad\text{and}\quad\sigma(D_n)\subset\mathcal{G}_\infty\quad\text{for
all }n,\text{ where }\mathcal{G}_\infty=\sigma(\cup_n\mathcal{G}_n).
\end{equation*}
If $C_n\rightarrow M$ stably and $D_n\rightarrow N$ stably in strong
sense, with respect to $(\mathcal{G}_n)$, then
\begin{equation*}
(C_n,D_n)\longrightarrow M\times N\quad\text{stably}.
\end{equation*} (Here, $M\times N$ is the kernel on $S\times S$ such
that $\bigl(M\times N\bigr)(\omega)=M(\omega)\times N(\omega)$
for all $\omega$).
\end{lem}
\begin{proof}
By standard arguments, since $S$ is separable and
$\sigma(C_n,D_n)\subset\mathcal{G}_\infty$, it suffices to prove
that $E\bigl\{I_H\, f_1(C_n)\,f_2(D_n)\}\rightarrow E\bigl\{I_H\,
M(f_1)\,N(f_2)\}$ whenever $H\in\cup_n\mathcal{G}_n$ and
$f_1,\,f_2\in C_b(S)$. Let $L_n=
E\bigl(f_2(D_n)\mid\mathcal{G}_n\bigr)-N(f_2)$. Since
$H\in\cup_n\mathcal{G}_n$, there is $k$ such that
$H\in\mathcal{G}_n$ for $n\geq k$. Thus,
\begin{gather*}
E\bigl\{I_H\, f_1(C_n)\,f_2(D_n)\}=E\bigl\{I_H\,
f_1(C_n)\,E\bigl(f_2(D_n)\mid\mathcal{G}_n\bigr)\}
\\=E\bigl\{I_H\,f_1(C_n)\,N(f_2)\}+E\bigl\{I_H\, f_1(C_n)\,L_n\}\quad\text{for all }n\geq k.
\end{gather*}
Finally, $\abs{E\bigl\{I_H\,
f_1(C_n)\,L_n\}\,}\leq\sup{\abs{f_1}}\,E\abs{L_n}\rightarrow 0$,
since $D_n\rightarrow N$ stably in strong sense, and
$E\bigl\{I_H\,f_1(C_n)\,N(f_2)\}\rightarrow
E\bigl\{I_H\,M(f_1)\,N(f_2)\}$ as $C_n\rightarrow M$ stably.
\end{proof}

\section{Main result}\label{main}

In the sequel, $(X_n:n\geq 1)$ is a sequence of real random
variables on the probability space $(\Omega,\mathcal{A},P)$ and
$(\mathcal{G}_n:n\geq 0)$ an (increasing) filtration. We assume
$E\abs{X_n}<\infty$ and we let
\begin{equation*}
Z_n=E(X_{n+1}\mid\mathcal{G}_n)\quad\text{and}\quad
\overline{X}_n=\frac{1}{n} \sum_{k=1}^nX_k.
\end{equation*}

In case $\sup_nEX_n^2<\infty$ and
\begin{equation}\label{basic}
E\bigr\{\bigl(E(Z_{n+1}\mid\mathcal{G}_n)-Z_n\bigr)^2\bigr\}=o(n^{-3}),
\end{equation}
the sequence $(Z_n)$ is an uniformly integrable quasi-martingale;
see e.g. page 532 of \cite{K}. Accordingly,
\begin{equation*}
Z_n\overset{a.s.,L_1}\longrightarrow Z
\end{equation*}
for some real random variable $Z$. Define
\begin{gather*}
C_n=\sqrt{n}\,\bigl(\overline{X}_n-Z_n\bigr),\quad
D_n=\sqrt{n}\,\bigl(Z_n-Z\bigr).
\end{gather*}

Let $\mathcal{N}(a,b)$ denote the one-dimensional Gaussian law with
mean $a$ and variance $b\geq 0$ (where $\mathcal{N}(a,0)=\delta_a$).
Note that $\mathcal{N}(0,L)$ is a kernel on $\mathbb{R}$ for each
real non negative random variable $L$. We are now in a position to
state our main result.

\begin{thm}\label{main}
Suppose $\sigma(X_n)\subset\mathcal{G}_n$ for each $n\geq 1$,
$(X_n^2)$ is uniformly integrable and condition \eqref{basic} holds.
Let us consider the following conditions

\vspace{0.2cm}

\begin{itemize}

\item[{\bf (a)}]
$\frac{1}{\sqrt{n}}\,E\bigl\{\max_{1\leq k\leq
n}k\,\abs{Z_{k-1}-Z_k}\bigr\}\longrightarrow 0$,

\vspace{0.2cm}

\item[{\bf (b)}] $\frac{1}{n}\sum_{k=1}^n\bigl\{X_k-Z_{k-1}+k(Z_{k-1}-Z_k)\bigr\}^2\overset{P}\longrightarrow
U$,

\vspace{0.2cm}

\item[{\bf (c)}] $\sqrt{n}\,E\bigl\{\sup_{k\geq n}\abs{Z_{k-1}-Z_k}\,\bigr\}\longrightarrow
0$,

\vspace{0.2cm}

\item[{\bf (d)}] $n\sum_{k\geq n}(Z_{k-1}-Z_k)^2\overset{P}\longrightarrow
V$,

\end{itemize}

\vspace{0.2cm}

\noindent where $U$ and $V$ are real non negative random variables.
Then, $C_n\rightarrow\mathcal{N}(0,U)$ stably under (a)-(b), and
$D_n\rightarrow\mathcal{N}(0,V)$ stably in strong sense, with
respect to $(\mathcal{G}_n)$, under (c)-(d). In
particular,\begin{equation*} (C_n,D_n)\longrightarrow
\mathcal{N}(0,U)\times\mathcal{N}(0,V)\quad\text{stably under
(a)-(b)-(c)-(d)}.
\end{equation*}
\end{thm}

\begin{proof} Since $\sigma(C_n)\subset\mathcal{G}_n$ and $Z$ can be taken
$\mathcal{G}_\infty$-measurable, Lemma \ref{plm} applies. Thus, it
suffices to prove that $C_n\rightarrow\mathcal{N}(0,U)$ stably and
$D_n\rightarrow\mathcal{N}(0,V)$ stably in strong sense.

\vspace{0.2cm}

{\bf "$C_n\rightarrow\mathcal{N}(0,U)$ stably".} Suppose conditions
(a)-(b) hold. First note that
\begin{gather*}
\sqrt{n}\,C_n=n\,\overline{X}_n-n\,Z_n=\sum_{k=1}^nX_k+\sum_{k=1}^n\bigl((k-1)Z_{k-1}-kZ_k\bigr)
\\=\sum_{k=1}^n\bigl\{X_k-Z_{k-1}+k(Z_{k-1}-Z_k)\bigr\}.
\end{gather*}
Letting
\begin{equation*}
Y_{n,k}=\frac{X_k-Z_{k-1}+k\bigl(E(Z_k\mid\mathcal{G}_{k-1})-Z_k\bigr)}{\sqrt{n}}\quad\text{and}\quad
Q_n=\frac{1}{\sqrt{n}}\,\sum_{k=1}^nk\bigl(Z_{k-1}-E(Z_k\mid\mathcal{G}_{k-1})\bigr),
\end{equation*}
it follows that $C_n=\sum_{k=1}^nY_{n,k}\,+\,Q_n$. By
\eqref{basic},
\begin{gather*}
E\abs{Q_n}\leq\frac{1}{\sqrt{n}}\,\sum_{k=1}^nk\,\sqrt{E\bigr\{\bigl(Z_{k-1}-E(Z_k\mid\mathcal{G}_{k-1})\bigr)^2\bigr\}}
=\frac{1}{\sqrt{n}}\,\sum_{k=1}^n\text{o}(k^{-1/2})\longrightarrow
0.
\end{gather*}
Hence, it suffices to prove that
$\sum_{k=1}^nY_{n,k}\rightarrow\mathcal{N}(0,U)$ stably. Letting
$\mathcal{F}_{n,k}=\mathcal{G}_k$, $k=1,\ldots,n$, one obtains $E\bigl(Y_{n,k}\mid\mathcal{F}_{n,k-1}\bigr)=0$ a.s.. Thus, by Corollary 7 of \cite{CLP}, $\sum_{k=1}^nY_{n,k}\rightarrow\mathcal{N}(0,U)$ stably whenever
\begin{equation*}
\text{(i)}\,\,\,E\bigl\{\max_{1\leq k\leq
n}\abs{Y_{n,k}}\bigr\}\longrightarrow 0;\quad
\text{(ii)}\,\,\,\sum_{k=1}^nY_{n,k}^2\overset{P}\longrightarrow U.
\end{equation*}

As to (i), first note that
\begin{equation*}
\sqrt{n}\max_{1\leq k\leq n}\abs{Y_{n,k}}\leq\max_{1\leq k\leq
n}\abs{X_k-Z_{k-1}}\,+\,\sum_{k=1}^nk\,\abs{E(Z_k\mid\mathcal{G}_{k-1})-Z_{k-1}}\,+\,
\max_{1\leq k\leq n}k\,\abs{Z_{k-1}-Z_k}.
\end{equation*}
Since $(X_n^2)$ is uniformly integrable, $((X_n-Z_{n-1})^2)$ is
uniformly integrable as well, and this implies
$\frac{1}{n}\,E\bigl\{\max_{1\leq k\leq
n}(X_k-Z_{k-1})^2\bigr\}\longrightarrow 0$. By condition
\eqref{basic},
\begin{equation*}
\frac{1}{\sqrt{n}}\sum_{k=1}^nk\,E\Abs{E(Z_k\mid\mathcal{G}_{k-1})-Z_{k-1}}=
\frac{1}{\sqrt{n}}\,\sum_{k=1}^n\text{o}(k^{-1/2})\longrightarrow 0.
\end{equation*}
Thus, (i) follows from condition (a).

As to (ii), write
\begin{gather*}
\sum_{k=1}^nY_{n,k}^2=\frac{1}{n}\sum_{k=1}^n\bigl(X_k-Z_{k-1}+k(Z_{k-1}-Z_k)\bigr)^2+\,
\frac{1}{n}\sum_{k=1}^nk^2\bigl(E(Z_k\mid\mathcal{G}_{k-1})-Z_{k-1}\bigr)^2+
\\+\,\frac{2}{n}\sum_{k=1}^n\bigl(X_k-Z_{k-1}+k(Z_{k-1}-Z_k)\bigr)\,k\,\bigl(E(Z_k\mid\mathcal{G}_{k-1})-Z_{k-1}\bigr)
\\=R_n+S_n+T_n\quad\text{say}.
\end{gather*}
Then, $R_n\overset{P}\rightarrow U$ by (b) and
$E\abs{S_n}=ES_n\rightarrow 0$ by \eqref{basic}. Further
$T_n\overset{P}\longrightarrow 0$, since
\begin{gather*}
\frac{T_n^2}{4}\leq
\frac{1}{n}\sum_{k=1}^n\bigl(X_k-Z_{k-1}+k(Z_{k-1}-Z_k)\bigr)^2\,\cdot\,\frac{1}{n}\sum_{k=1}^nk^2\bigl(E(Z_k\mid\mathcal{G}_{k-1})-Z_{k-1}\bigr)^2=R_n\,S_n.
\end{gather*}
Hence, (ii) holds, and this concludes the proof of $C_n\rightarrow\mathcal{N}(0,U)$ stably.

\vspace{0.4cm}

{\bf "$D_n\rightarrow\mathcal{N}(0,V)$ stably in strong sense".}
Suppose conditions (c)-(d) hold. We first recall a known result; see
Example 6 of \cite{CLP}. Let $(L_n)$ be a
$(\mathcal{G}_n)$-martingale such that
$L_n\overset{a.s.,L_1}\longrightarrow L$ for some real random
variable $L$. Then,
\begin{gather*}
\sqrt{n}\,\bigl(L_n-L\bigr)\longrightarrow\mathcal{N}(0,V)\quad\text{stably
in strong sense with respect to }(\mathcal{G}_n),
\end{gather*}
provided
\begin{equation*}
\text{(c*)}\,\,\,\sqrt{n}\,E\bigl\{\sup_{k\geq
n}\abs{L_{k-1}-L_k}\,\bigr\}\longrightarrow 0;\quad\text{(d*)}\,\,\,
n\sum_{k\geq n}(L_{k-1}-L_k)^2\overset{P}\longrightarrow V.
\end{equation*}
Next, define $L_0=Z_0$ and
\begin{gather*}
L_n=Z_n-\sum_{k=0}^{n-1}\bigl(E(Z_{k+1}\mid\mathcal{G}_k)-Z_k\bigr).
\end{gather*}
Then, $(L_n)$ is a $(\mathcal{G}_n)$-martingale. Also,
$L_n\overset{a.s.,L_1}\longrightarrow L$ for some $L$, as $(Z_n)$ is
an uniformly integrable quasi martingale. In particular, $L_n-L$ can
be written as $L_n-L=\sum_{k\geq n}(L_k-L_{k+1})$ a.s.. Similarly,
$Z_n-Z=\sum_{k\geq n}(Z_k-Z_{k+1})$ a.s.. It follows that
\begin{gather*}
E\Abs{D_n-\sqrt{n}(L_n-L)}=\sqrt{n}\,E\Abs{(Z_n-Z)-(L_n-L)}
\\=\sqrt{n}\,E\Abs{\,\sum_{k\geq
n}\bigl\{(Z_k-L_k)-(Z_{k+1}-L_{k+1})\bigr\}}
\\\leq\sqrt{n}\,\sum_{k\geq
n}E\Abs{Z_k-E(Z_{k+1}\mid\mathcal{G}_k)}=\sqrt{n}\,\sum_{k\geq
n}\text{o}(k^{-3/2})\longrightarrow 0.
\end{gather*}
Thus, $D_n\rightarrow\mathcal{N}(0,V)$ stably in strong sense if and
only if $\sqrt{n}(L_n-L)\rightarrow\mathcal{N}(0,V)$ stably in
strong sense, and to conclude the proof it suffices to check
conditions (c*)-(d*). In turn, (c*)-(d*) are a straightforward
consequence of conditions \eqref{basic}, (c), (d) and
\begin{equation*}
L_{k-1}-L_k=\bigl(Z_{k-1}-Z_k\bigr)+\bigl(E(Z_k\mid\mathcal{G}_{k-1})-Z_{k-1}\bigr).
\end{equation*}
\end{proof}

Some remarks on Theorem \ref{main} are in order.

In real problems, one of the quantities of main interest is
\begin{equation*}
W_n=\sqrt{n}\,\bigl(\overline{X}_n-Z).
\end{equation*}
And, under the assumptions of Theorem \ref{main}, one obtains
\begin{equation*}
W_n=C_n+D_n\longrightarrow\mathcal{N}(0,U+V)\quad\text{stably}.
\end{equation*}

Condition \eqref{basic} trivially holds when $(X_n)$ is
conditionally identically distributed, in the sense of \cite{BPR},
with respect to the filtration $(\mathcal{G}_n)$. In this case, in
fact, $(Z_n)$ is even a $(\mathcal{G}_n)$-martingale. In particular,
\eqref{basic} holds if $(X_n)$ is exchangeable and
$\mathcal{G}_n=\sigma(X_1,\ldots,X_n)$.

Under (c), condition (a) can be replaced by

\vspace{0.2cm}

\begin{itemize}
\item[{\bf (a*)}]
$\sup_n\,\frac{1}{n}\sum_{k=1}^nk^2E\bigl\{(Z_{k-1}-Z_k)^2\bigr\}<\infty$.
\end{itemize}

\vspace{0.2cm}

\noindent Indeed, (a*) and (c) imply (a) (we omit calculations). Note that, for proving $C_n\rightarrow\mathcal{N}(0,U)$ stably under (a*)-(b)-(c), one can rely on more classical versions of the martingale CLT, such as Theorem 3.2 of \cite{HH}.

To check conditions (b) and (d), the following simple lemma can
help.

\begin{lem}\label{hftuimn}
Let $(Y_n)$ be a $(\mathcal{G}_n)$-adapted sequence of real random
variables. If $\sum_{n=1}^\infty\frac{EY_n^2}{n^2}<\infty$ and
$E\bigl(Y_{n+1}\mid\mathcal{G}_n\bigr)\overset{a.s.}\longrightarrow
Y$, for some random variable $Y$, then
\begin{equation*}
n\sum_{k\geq n}\frac{Y_k}{k^2}\,\overset{a.s.}\longrightarrow
Y\quad\text{and}\quad\frac{1}{n}\sum_{k=1}^nY_k\overset{a.s.}\longrightarrow
Y.
\end{equation*}
\end{lem}
\begin{proof}
Let
$L_n=\sum_{k=1}^n\frac{Y_k-E\bigl(Y_k\mid\mathcal{G}_{k-1}\bigr)}{k}$.
Then, $L_n$ is a $(\mathcal{G}_n)$-martingale such that
\begin{equation*}
\sup_nEL_n^2\leq 4\,\sum_k\frac{EY_k^2}{k^2}<\infty.
\end{equation*} Thus, $L_n$ converges a.s. and Abel summation formula yields
\begin{equation*}
n\sum_{k\geq
n}\,\frac{Y_k-E\bigl(Y_k\mid\mathcal{G}_{k-1}\bigr)}{k^2}\,\overset{a.s.}\longrightarrow
0.
\end{equation*}
Since
$E\bigl(Y_{n+1}\mid\mathcal{G}_n\bigr)\overset{a.s.}\longrightarrow
Y$ and $n\sum_{k\geq n}\,\frac{1}{k^2}\longrightarrow 1$, it follows
that
\begin{gather*}
n\sum_{k\geq n}\,\frac{Y_k}{k^2}=n\sum_{k\geq
n}\,\frac{Y_k-E\bigl(Y_k\mid\mathcal{G}_{k-1}\bigr)}{k^2}\,+\,n\sum_{k\geq
n}\,\frac{E\bigl(Y_k\mid\mathcal{G}_{k-1}\bigr)}{k^2}\overset{a.s.}\longrightarrow
Y.
\end{gather*}
Similarly, Kroneker lemma and
$E\bigl(Y_{n+1}\mid\mathcal{G}_n\bigr)\overset{a.s.}\longrightarrow
Y$ yield
\begin{gather*}
\frac{1}{n}\sum_{k=1}^nY_k=\frac{1}{n}\sum_{k=1}^nE(Y_k\mid\mathcal{G}_{k-1})\,+\,\frac{1}{n}\sum_{k=1}^nk\,\frac{Y_k-E\bigl(Y_k\mid\mathcal{G}_{k-1}\bigr)}{k}\overset{a.s.}\longrightarrow
Y.
\end{gather*}
\end{proof}

Our last comment needs a formal remark.

\begin{rem}\label{stabqc}
As regards $D_n$, a natural question is whether
\begin{equation}\label{pevb}
E\bigl(f(D_n)\mid\mathcal{G}_n\bigr)\overset{a.s.}\longrightarrow
\mathcal{N}(0,V)(f)\quad\text{for each }f\in C_b(\mathbb{R}).
\end{equation}
This is a strengthening of $D_n\rightarrow\mathcal{N}(0,V)$ stably
in strong sense, as $E\bigl(f(D_n)\mid\mathcal{G}_n\bigr)$ is
requested to converge a.s. and not only in probability. Let $(X_n)$
be a (non necessarily $(\mathcal{G}_n)$-adapted) sequence of
integrable random variables. Then, for \eqref{pevb} to be true, it
is enough that $(Z_n)$ is uniformly integrable and
\begin{gather*}
\sum_{k\geq
1}\sqrt{k}\,E\Abs{\,E(Z_k\mid\mathcal{G}_{k-1})-Z_{k-1}}<\infty,
\\E\bigl\{\sup_{k\geq
1}\sqrt{k}\,\abs{Z_{k-1}-Z_k}\bigr\}<\infty,\quad n\sum_{k\geq
n}(Z_{k-1}-Z_k)^2\overset{a.s.}\longrightarrow V.
\end{gather*}
The proof is essentially the same as that of Theorem \ref{main}, up
to using Theorem 2.2 of \cite{C} instead of Example 6 of \cite{CLP}.
\end{rem}

\section{Applications}\label{appl}

This section is split into four subsections, arranged in increasing
order of length.

\subsection{$r$-step predictions}\label{kstep} Suppose we are requested to make conditional forecasts on a sequence of events $A_n\in\mathcal{G}_n$. To fix ideas, for each $n$, we aim to predict
\begin{equation*}
A_n^*=\bigl(\cap_{j\in J}A_{n+j}\bigr)\cap\bigl(\cap_{j\in J^c}A_{n+j}^c\bigr)
\end{equation*}
conditionally on $\mathcal{G}_n$, where $J$ is a given subset of
$\{1,\ldots,r\}$ and $J^c=\{1,\ldots,r\}\setminus J$. Letting
$X_n=I_{A_n}$, the predictor can be written as
\begin{gather*}
Z_n^*=E\bigl\{\,\prod_{j\in J}X_{n+j}\,\prod_{j\in J^c}(1-X_{n+j})\mid\mathcal{G}_n\bigr\}.
\end{gather*}

In the spirit of Section \ref{intro}, when $Z_n^*$ cannot be
evaluated in closed form, one needs to estimate it. Under some
assumptions, in particular when $(X_n)$ is exchangeable and
$\mathcal{G}_n=\sigma(X_1,\ldots,X_n)$, a reasonable estimate of
$Z_n^*$ is $\overline{X}_n^h(1-\overline{X}_n)^{r-h}$ where
$h=\,$card$(J)$. Usually, under such assumptions, one also has
$Z_n\overset{a.s.}\longrightarrow Z$ and
$Z_n^*\overset{a.s.}\longrightarrow Z^h(1-Z)^{r-h}$ for some random
variable $Z$. So, it makes sense to define
\begin{gather*}
C_n^*=\sqrt{n}\,\bigl\{\overline{X}_n^h(1-\overline{X}_n)^{r-h}-Z_n^*\bigr\},\quad
D_n^*=\sqrt{n}\,\bigl\{Z_n^*-Z^h(1-Z)^{r-h}\bigr\}.
\end{gather*}

Next result is a straightforward consequence of Theorem \ref{main}.

\begin{cor}
Let $(X_n)$ be a $(\mathcal{G}_n)$-adapted sequence of indicators
satisfying \eqref{basic}. If conditions (a)-(b)-(c)-(d) of Theorem
\ref{main} hold, then \begin{gather*} (C_n^*,D_n^*)\longrightarrow
\mathcal{N}(0,\sigma^2U)\times\mathcal{N}(0,\sigma^2V)\quad\text{stably,
where}
\\\sigma^2=\bigl\{h\,Z^{h-1}(1-Z)^{r-h}-(r-h)\,Z^h(1-Z)^{r-h-1}\bigr\}^2.
\end{gather*}
\end{cor}
\begin{proof}
We just give a sketch of the proof. Let $f(x)=x^h(1-x)^{r-h}$.
Basing on (c), it can be shown that
$\sqrt{n}\,E\Abs{\,Z_n^*-f(Z_n)}\longrightarrow 0$. Thus, $C_n^*$
can be replaced by
$\sqrt{n}\,\bigl\{f(\overline{X}_n)-f(Z_n)\bigr\}$ and $D_n^*$ by
$\sqrt{n}\,\bigl\{f(Z_n)-f(Z)\bigr\}$. By the mean value theorem,
\begin{equation*}
\sqrt{n}\,\bigl\{f(\overline{X}_n)-f(Z_n)\bigr\}=\sqrt{n}\,f'(M_n)\,(\overline{X}_n-Z_n)=f'(M_n)\,C_n
\end{equation*}
where $M_n$ is between $\overline{X}_n$ and $Z_n$. By \eqref{basic},
$Z_n\overset{a.s}\longrightarrow Z$ and
$\overline{X}_n\overset{a.s}\longrightarrow Z$. Hence,
$f'(M_n)\overset{a.s}\longrightarrow f'(Z)$ as $f'$ is continuous.
By Theorem \ref{main}, $C_n\rightarrow \mathcal{N}(0,U)$ stably.
Thus,
\begin{equation*}
\sqrt{n}\,\bigl\{f(\overline{X}_n)-f(Z_n)\bigr\}\longrightarrow
f'(Z)\,\mathcal{N}(0,U)=\mathcal{N}(0,\sigma^2U)\quad\text{stably}.
\end{equation*}
By a similar argument, it can be seen that
$\sqrt{n}\,\bigl\{f(Z_n)-f(Z)\bigr\}\longrightarrow\mathcal{N}(0,\sigma^2V)$
stably in strong sense. An application of Lemma \ref{plm} concludes
the proof.
\end{proof}

\subsection{Poisson-Dirichlet sequences}\label{poidir}
Let $\mathcal{Y}$ be a finite set and $(Y_n)$ a sequence of
$\mathcal{Y}$-valued random variables satisfying
\begin{equation*}
P\bigl(Y_{n+1}\in A\mid Y_1,\ldots,Y_n\bigr)=\frac{\sum_{y\in
A}(S_{n,y}-\alpha)\,I_{\{S_{n,y}\neq
0\}}+\bigl(\theta+\alpha\sum_{y\in\mathcal{Y}}I_{\{S_{n,y}\neq
0\}}\bigr)\,\nu(A)}{\theta+n}
\end{equation*}
a.s. for all $A\subset\mathcal{Y}$ and $n\geq 1$. Here, $0\leq
\alpha<1$ and $\theta>-\alpha$ are constants, $\nu$ is the
probability distribution of $Y_1$ and
$S_{n,y}=\sum_{k=1}^nI_{\{Y_k=y\}}$.

Sequences $(Y_n)$ of this type play a role in various frameworks, mainly in population-genetics. They can be regarded as a generalization of those exchangeable sequences directed by a two parameter Poisson-Dirichlet process; see \cite{PY}. For $\alpha=0$, $(Y_n)$ reduces to a classical Dirichlet sequence (i.e., an exchangeable sequence directed by a Dirichlet process). But, for $\alpha\neq 0$, $(Y_n)$ may even fail to be exchangeable.

From the point of view of Theorem \ref{main}, however, the only important thing is that $P\bigl(Y_{n+1}\in\cdot\mid Y_1,\ldots,Y_n\bigr)$ can be written down explicitly. Indeed, the following result is available.

\begin{cor} Let $\mathcal{G}_n=\sigma(Y_1,\ldots,Y_n)$ and $X_n=I_A(Y_n)$, where $A\subset\mathcal{Y}$.
Then, condition \eqref{basic} holds (so that
$Z_n\overset{a.s.}\longrightarrow Z$) and
\begin{equation*}
(C_n,D_n)\longrightarrow\delta_0\times\mathcal{N}\bigl(0,Z(1-Z)\bigr)\quad\text{stably}.
\end{equation*}
\end{cor}
\begin{proof}
Let $Q_n=-\alpha\,\sum_{y\in A}I_{\{S_{n,y}\neq 0\}}+\bigl(\theta+\alpha\sum_{y\in\mathcal{Y}}I_{\{S_{n,y}\neq 0\}}\bigr)\,\nu(A)$. Since
\begin{equation*}
Z_n=P\bigl(Y_{n+1}\in A\mid Y_1,\ldots,Y_n\bigr)=\frac{n\,\overline{X}_n\,+\,Q_n}{\theta+n}\quad\text{and}\quad \abs{Q_n}\leq c
\end{equation*}
for some constant $c$, then $C_n\overset{a.s.}\longrightarrow 0$. By Lemma \ref{plm} and Theorem \ref{main}, thus, it suffices to check conditions \eqref{basic}, (c) and (d) with $V=Z(1-Z)$. On noting that
\begin{equation*}
Z_{n+1}-Z_n=\frac{X_{n+1}-Z_n}{\theta+n+1}\,+\,\frac{Q_{n+1}-Q_n}{\theta+n+1},
\end{equation*}
condition (c) trivially holds. Since $S_{n+1,y}=S_{n,y}+I_{\{Y_{n+1}=y\}}$, one obtains
\begin{gather*}
Q_{n+1}-Q_n=-\alpha\,\nu(A^c)\,\sum_{y\in A}I_{\{S_{n,y}=0\}}I_{\{Y_{n+1}=y\}}\,+\,\alpha\,\nu(A)\,\sum_{y\in A^c}I_{\{S_{n,y}=0\}}I_{\{Y_{n+1}=y\}}.
\end{gather*}
It follows that
\begin{equation*}
E\bigl\{\abs{Q_{n+1}-Q_n}\mid\mathcal{G}_n\bigr\}\leq
2\,\sum_{y\in\mathcal{Y}}I_{\{S_{n,y}=0\}}\,P\bigl(Y_{n+1}=y\mid\mathcal{G}_n\bigr)\leq\frac{d}{\theta+n}\quad\text{a.s.}
\end{equation*}
for some constant $d$, and this implies
\begin{equation*}
\Abs{E\bigl(Z_{n+1}\mid\mathcal{G}_n\bigr)-Z_n}=\frac{\Abs{E\bigl(Q_{n+1}-Q_n\mid\mathcal{G}_n\bigr)}}{\theta+n+1}\leq\frac{d}{(\theta+n)^2}\quad\text{a.s.}.
\end{equation*}
Hence, condition \eqref{basic} holds. To check (d), note that
$\sum_kk^2E\bigl\{(Z_{k-1}-Z_k)^4\bigr\}<\infty$. Since
$Z_k\overset{a.s.}\longrightarrow Z$ (by \eqref{basic}) one also
obtains
\begin{gather*}
E\bigl\{(X_k-Z_{k-1})^2\mid\mathcal{G}_{k-1}\bigr\}=Z_{k-1}-Z_{k-1}^2\overset{a.s.}\longrightarrow
Z(1-Z),
\\E\bigl\{(Q_k-Q_{k-1})^2\mid\mathcal{G}_{k-1}\bigr\}\,+\,
2\,E\bigl\{(X_k-Z_{k-1})\,(Q_k-Q_{k-1})\mid\mathcal{G}_{k-1}\bigr\}\overset{a.s.}\longrightarrow
0.
\end{gather*}
Thus,
$k^2E\bigl\{(Z_{k-1}-Z_k)^2\mid\mathcal{G}_{k-1}\bigr\}\overset{a.s.}\longrightarrow
Z(1-Z)$. Letting $Y_k=k^2(Z_{k-1}-Z_k)^2$ and $Y=Z(1-Z)$, Lemma
\ref{hftuimn} implies \begin{equation*} n\sum_{k\geq
n}(Z_{k-1}-Z_k)^2=n\sum_{k\geq
n}\frac{Y_k}{k^2}\overset{a.s.}\longrightarrow Z(1-Z).
\end{equation*}
\end{proof}

As it is clear from the previous proof, all conditions of Remark
\ref{stabqc} are satisfied. Therefore, $D_n$ meets condition
\eqref{pevb} with $V=Z(1-Z)$.

\subsection{Two color randomly reinforced generalized Polya urns}\label{2col}
An urn contains $b>0$ black balls and $r>0$ red balls. At each time
$n\geq 1$, a ball is drawn and then replaced together with a random
number of balls of the same color. Say that $B_n$ black balls or
$R_n$ red balls are added to the urn according to whether $X_n=1$ or
$X_n=0$, where $X_n$ is the indicator of $\{$black ball at time
$n\}$.

Urns of this type have some history: see \cite{AMS}, \cite{BCL}, \cite{BPR},
\cite{C}, \cite{MF}, \cite{P} and references therein.

To model such urns, we assume $X_n,\,B_n,\,R_n$ random variables on
the probability space $(\Omega,\mathcal{A},P)$ such that

\vspace{0.2cm}

\begin{itemize}
\item[{\bf ($*$)}] $X_n\in\{0,1\}$, $\,\,\,\,B_n\geq 0$, $\,\,\,\,R_n\geq 0$,
\begin{gather*}
(B_n,R_n)\text{ independent of }\,
\bigl(X_1,B_1,R_1,\ldots,X_{n-1},B_{n-1},R_{n-1},X_n\bigr),
\\Z_n=P\bigl(X_{n+1}=1\mid\mathcal{G}_n\bigr)=\frac{b+\sum_{k=1}^nB_kX_k}{b+r+\sum_{k=1}^n\bigl(B_kX_k+R_k(1-X_k)\bigr)}\,\,\text{
a.s.},
\end{gather*}
for each $n\geq 1$, where
\begin{equation*}
\mathcal{G}_0=\{\emptyset,\Omega\},\quad\mathcal{G}_n=\sigma\bigl(X_1,B_1,R_1,\ldots,X_n,B_n,R_n\bigr).
\end{equation*}
\end{itemize}

\vspace{0.2cm}

In the particular case $B_n=R_n$, in Example 3.5 of \cite{BPR}, it
is shown that $C_n$ converges stably to a Gaussian kernel whenever
$EB_1^2<\infty$ and $B_n\sim B_1$ for all $n$. Further, in Corollary
4.1 of \cite{C}, $D_n$ is shown to satisfy condition \eqref{pevb}.
The latter result on $D_n$ is extended to $B_n\neq R_n$ in
\cite{AMS}, under the assumptions that $B_1+R_1$ has compact
support, $EB_1=ER_1$, and $(B_n,R_n)\sim (B_1,R_1)$ for all $n$.

Basing on Theorem \ref{main}, condition \eqref{pevb} can be shown to
hold more generally. Indeed, it is fundamental that $EB_n=ER_n$ for
all $n$ and the three sequences $(EB_n)$, $(EB_n^2)$, $(ER_n^2)$
approach a limit. But identity in distribution of $(B_n,R_n)$ can be
dropped and compact support of $B_n+R_n$ can be replaced by a moment
condition such as
\begin{equation}\label{momcon}
\sup_n E\bigl\{(B_n+R_n)^u\bigr\}<\infty\quad\text{for some }u>2.
\end{equation}
Under these conditions, not only $D_n$ meets \eqref{pevb}, but the
pairs $(C_n,D_n)$ converge stably as well. In particular, one
obtains stable convergence of $W_n=C_n+D_n$ which is of potential
interest in urn problems.

\begin{cor}\label{poi}
In addition to ($*$) and \eqref{momcon}, suppose $EB_n=ER_n$ for all
$n$ and
\begin{gather*}
m:=\lim_nEB_n>0,\quad q:=\lim_nEB_n^2,\quad s:=\lim_nER_n^2.
\end{gather*}
Then, condition \eqref{basic} holds (so that
$Z_n\overset{a.s.}\longrightarrow Z$) and
\begin{gather*}
(C_n,D_n)\longrightarrow
\mathcal{N}(0,U)\times\mathcal{N}(0,V)\quad\text{stably, where}
\\U=Z(1-Z)\,\bigl(\frac{(1-Z)q+Zs}{m^2}-1\bigr)\quad\text{and}\quad
V=Z(1-Z)\,\frac{(1-Z)q+Zs}{m^2}.
\end{gather*}
In particular, $W_n=C_n+D_n\longrightarrow\mathcal{N}(0,U+V)$
stably. Moreover, $D_n$ meets condition \eqref{pevb}, that is,
$E\bigl(f(D_n)\mid\mathcal{G}_n\bigr)\overset{a.s.}\longrightarrow
\mathcal{N}(0,V)(f)$ for each $f\in C_b(\mathbb{R})$.
\end{cor}

It is worth noting that, arguing as in \cite{AMS} and \cite{MF}, one
obtains $P(Z=z)=0$ for all $z$. Thus, $\mathcal{N}(0,V)$ is a non
degenerate kernel. In turn, $\mathcal{N}(0,U)$ is non degenerate
unless $q=s=m^2$, and this happens if and only if both $B_n$ and
$R_n$ converge in probability (necessarily to $m$). In the latter
case ($q=s=m^2$), $C_n\overset{P}\longrightarrow 0$ and condition
\eqref{pevb} holds with $V=Z(1-Z)$. Thus, in a sense, randomly
reinforced urns behave as classical Polya urns (i.e., those urns
with $B_n=R_n=m$) whenever the reinforcements converge in
probability.

The proof of Corollary \ref{poi} is deferred to the Appendix as it
needs some work. Here, to point out the underlying argument, we
sketch such a proof under the superfluous but simplifying assumption
that $B_n\vee R_n\leq c$ for all $n$ and some constant $c$. Let
\begin{equation*}
S_n=b+r+\sum_{k=1}^n\bigl(B_kX_k+R_k(1-X_k)\bigr).
\end{equation*}

After some algebra, $Z_{n+1}-Z_n$ can be written as
\begin{gather*}
Z_{n+1}-Z_n=\frac{(1-Z_n)\,X_{n+1}\,B_{n+1}\,-\,Z_n\,(1-X_{n+1})\,R_{n+1}}{S_{n+1}}
\\=\frac{(1-Z_n)\,X_{n+1}\,B_{n+1}}{S_n+B_{n+1}}\,-\,\frac{Z_n\,(1-X_{n+1})\,R_{n+1}}{S_n+R_{n+1}}.
\end{gather*}
By ($*$) and $EB_{n+1}=ER_{n+1}$,
\begin{gather*}
E\bigl(Z_{n+1}-Z_n\mid\mathcal{G}_n\bigr)=Z_n(1-Z_n)\,E\bigl\{\,\frac{B_{n+1}}{S_n+B_{n+1}}-\frac{R_{n+1}}{S_n+R_{n+1}}
\mid\mathcal{G}_n\bigr\}
\\=Z_n(1-Z_n)\,E\bigl\{\,\frac{B_{n+1}}{S_n+B_{n+1}}-\frac{B_{n+1}}{S_n}-\frac{R_{n+1}}{S_n+R_{n+1}}+\frac{R_{n+1}}{S_n}
\mid\mathcal{G}_n\bigr\}
\\=Z_n(1-Z_n)\,E\bigl\{\,-\frac{B_{n+1}^2}{S_n(S_n+B_{n+1})}\,+\,\frac{R_{n+1}^2}{S_n(S_n+R_{n+1})}\,\mid\mathcal{G}_n\bigr\}\quad\text{a.s.}.
\end{gather*}
Thus,
$\Abs{\,E\bigl(Z_{n+1}\mid\mathcal{G}_n\bigr)-Z_n}\leq\frac{EB_{n+1}^2+ER_{n+1}^2}{S_n^2}\,$
a.s.. Since $\sup_n\bigl(EB_n^2+ER_n^2\bigr)<\infty$ and
$E(S_n^{-p})=\,$O$(n^{-p})$ for all $p>0$ (as shown in Lemma
\ref{jrblk}) then
\begin{equation*}
E\bigl\{\abs{E(Z_{n+1}\mid\mathcal{G}_n)-Z_n}^p\bigr\}=\,\text{O}(n^{-2p})\quad\text{for
all }p>0.
\end{equation*}
In particular, condition \eqref{basic} holds and
$\sum_k\sqrt{k}\,E\Abs{\,E(Z_k\mid\mathcal{G}_{k-1})-Z_{k-1}}<\infty$.

To conclude the proof, in view of Lemma \ref{plm}, Theorem
\ref{main} and Remark \ref{stabqc}, it suffices to check conditions
(a), (b) and
\begin{equation*}
\text{(i)}\,\,\,E\bigl\{\sup_{k\geq
1}\sqrt{k}\,\abs{Z_{k-1}-Z_k}\bigr\}<\infty;\quad\text{(ii)}\,\,\,n\sum_{k\geq
n}(Z_{k-1}-Z_k)^2\overset{a.s.}\longrightarrow V.
\end{equation*}
Conditions (a) and (i) are straightforward consequences of
$\abs{Z_{n+1}-Z_n}\leq\frac{c}{S_n}$ and $E(S_n^{-p})=\,$O$(n^{-p})$
for all $p>0$. Condition (b) follows from the same argument as (ii).
And to prove (ii), it suffices to show that
$E(Y_{n+1}\mid\mathcal{G}_n)\overset{a.s.}\longrightarrow V$ where
$Y_n=n^2(Z_{n-1}-Z_n)^2$; see Lemma \ref{hftuimn}. Write
$(n+1)^{-2}E(Y_{n+1}\mid\mathcal{G}_n)$ as
\begin{gather*}
Z_n(1-Z_n)^2E\bigl\{\frac{B_{n+1}^2}{(S_n+B_{n+1})^2}\mid\mathcal{G}_n\bigr\}\,+\,Z_n^2(1-Z_n)E\bigl\{\frac{R_{n+1}^2}{(S_n+R_{n+1})^2}\mid\mathcal{G}_n\bigr\}.
\end{gather*}
Since $\frac{S_n}{n} \overset{a.s.}\longrightarrow m$ (by Lemma
\ref{jrblk}) and $B_{n+1}\leq c$, then
\begin{gather*}
n^2E\bigl\{\frac{B_{n+1}^2}{(S_n+B_{n+1})^2}\mid\mathcal{G}_n\bigr\}\leq n^2E\bigl\{\frac{B_{n+1}^2}{S_n^2}\mid\mathcal{G}_n\bigr\}=n^2\frac{EB_{n+1}^2}{S_n^2}\overset{a.s.}\longrightarrow\frac{q}{m^2}\,\text{ and}
\\n^2E\bigl\{\frac{B_{n+1}^2}{(S_n+B_{n+1})^2}\mid\mathcal{G}_n\bigr\}\geq n^2E\bigl\{\frac{B_{n+1}^2}{(S_n+c)^2}\mid\mathcal{G}_n\bigr\}=n^2\frac{EB_{n+1}^2}{(S_n+c)^2}\overset{a.s.}\longrightarrow\frac{q}{m^2}.
\end{gather*}
Similarly, $n^2E\bigl\{\frac{R_{n+1}^2}{(S_n+R_{n+1})^2}\mid\mathcal{G}_n\bigr\}\overset{a.s.}\longrightarrow\frac{s}{m^2}$. Since $Z_n\overset{a.s.}\longrightarrow Z$, it follows that
\begin{gather*}
E(Y_{n+1}\mid\mathcal{G}_n)\overset{a.s.}\longrightarrow \,\,Z(1-Z)^2\frac{q}{m^2}+Z^2(1-Z)\frac{s}{m^2}=V.
\end{gather*}
This concludes the (sketch of the) proof.

\begin{rem}
In order to $(C_n,D_n)\longrightarrow
\mathcal{N}(0,U)\times\mathcal{N}(0,V)$ stably, some of the
assumptions of Corollary \ref{poi} can be stated in a different
form. We mention two (independent) facts.

First, condition \eqref{momcon} can be weakened into uniform
integrability of $(B_n+R_n)^2$.

Second, $(B_n,R_n)$ independent of
$\mathcal{G}_{n-1}\vee\sigma(X_n)$ can be replaced by the following
four conditions:

\begin{itemize}
\item[(i)] $(B_n,R_n)$ conditionally independent of $X_n$ given
$\mathcal{G}_{n-1}$;
\item[(ii)] Condition \eqref{momcon} holds for some $u>4$;
\item[(iii)] There are an integer $n_0$ and a constant $l>0$ such that
\begin{gather*}
E\bigl(B_n\wedge n^{1/4}\mid\mathcal{G}_{n-1}\bigr)\geq
l\,\,\text{and}\,\,E\bigl(R_n\wedge
n^{1/4}\mid\mathcal{G}_{n-1}\bigr)\geq l\,\text{ a.s. whenever
}n\geq n_0; \end{gather*}
\item[(iv)] There are random variables $m,\,q,\,s$ such that
\begin{gather*}
E\bigl(B_n\mid\mathcal{G}_{n-1}\bigr)=E\bigl(R_n\mid\mathcal{G}_{n-1}\bigr)\overset{P}\longrightarrow
m,\quad
E\bigl(B_n^2\mid\mathcal{G}_{n-1}\bigr)\overset{P}\longrightarrow
q,\quad
E\bigl(R_n^2\mid\mathcal{G}_{n-1}\bigr)\overset{P}\longrightarrow s.
\end{gather*}

\end{itemize}

\noindent Even if in a different framework, conditions similar to
(i)-(iv) are in \cite{BH}.

\end{rem}

\subsection{The multicolor case}\label{multi}

To avoid technicalities, we firstly investigated two color urns, but
the results in Subsection \ref{2col} extend to the multicolor case.

An urn contains $a_j>0$ balls of color $j\in\{1,\ldots,d\}$ where
$d\geq 2$. Let $X_{n,j}$ denote the indicator of $\{$ball of color
$j$ at time $n\}$. In case $X_{n,j}=1$, the ball which has been
drawn is replaced together with $A_{n,j}$ more balls of color $j$.
Formally, we assume $\bigl\{X_{n,j},\,A_{n,j}:n\geq 1,\,1\leq j\leq
d\bigr\}$ random variables on the probability space
$(\Omega,\mathcal{A},P)$ satisfying

\vspace{0.2cm}

\begin{itemize}
\item[{\bf ($**$)}] $\,\,\,\,X_{n,j}\in\{0,1\}$,
$\,\,\,\,\sum_{j=1}^dX_{n,j}=1$, $\,\,\,\,A_{n,j}\geq 0$,
\begin{gather*}
(A_{n,1},\ldots,A_{n,d})\text{ independent of }\,
\bigl(A_{k,j},\,X_{k,j},\,X_{n,j}:1\leq k<n,\,1\leq j\leq d\bigr),
\\Z_{n,j}=P\bigl(X_{n+1,j}=1\mid\mathcal{G}_n\bigr)=\frac{a_j+\sum_{k=1}^nA_{k,j}X_{k,j}}{\sum_{i=1}^da_i+\sum_{k=1}^n\sum_{i=1}^dA_{k,i}X_{k,i}}\,\,\text{
a.s.},
\\\text{where}\quad\mathcal{G}_0=\{\emptyset,\Omega\},\quad\mathcal{G}_n=\sigma\bigl(A_{k,j},\,X_{k,j}:1\leq
k\leq n,\,1\leq j\leq d\bigr).
\end{gather*}
\end{itemize}

\vspace{0.2cm}

Note that
\begin{gather*}
Z_{n+1,j}-Z_{n,j}=(1-Z_{n,j})\,\frac{A_{n+1,j}\,X_{n+1,j}}{S_n+A_{n+1,j}}\,-\,Z_{n,j}\sum_{i\neq
j}\frac{A_{n+1,i}\,X_{n+1,i}}{S_n+A_{n+1,i}}
\\\text{where
}\,S_n=\sum_{i=1}^da_i+\sum_{k=1}^n\sum_{i=1}^dA_{k,i}X_{k,i}.
\end{gather*}

In addition to ($**$), as in Subsection \ref{2col}, we ask the
moment condition
\begin{equation}\label{hg5dr}
\sup_n
E\bigl\{\bigl(\,\sum_{j=1}^dA_{n,j}\bigr)^u\bigr\}<\infty\quad\text{for
some }u>2.
\end{equation}
Further, it is fundamental that
\begin{gather}\label{pawbn}
EA_{n,j}=EA_{n,1}\quad\text{for each }n\geq 1\text{ and }1\leq j\leq
d,\,\text{ and}
\\m:=\lim_nEA_{n,1}>0,\quad q_j:=\lim_nEA_{n,j}^2\quad\text{for each }1\leq j\leq d.\notag
\end{gather}

Fix $1\leq j\leq d$. Since $EA_{n,i}=EA_{n,1}$ for all $n$ and $i$,
the same calculation as in Subsection \ref{2col} yields
\begin{gather*}
\Abs{\,E\bigl(Z_{n+1,j}\mid\mathcal{G}_n\bigr)-Z_{n,j}}\leq\frac{\sum_{i=1}^dEA_{n+1,i}^2}{S_n^2}\quad\text{a.s.}.
\end{gather*}
Also, $E(S_n^{-p})=\,$O$(n^{-p})$ for all $p>0$; see Remark
\ref{chissa1}. Thus,
\begin{equation}\label{aw38j}
E\bigl\{\abs{E\bigl(Z_{n+1,j}\mid\mathcal{G}_n\bigr)-Z_{n,j}}^p\bigr\}=\,\text{O}(n^{-2p})\quad\text{for
all }p>0.
\end{equation}
In particular, $Z_{n,j}$ meets condition \eqref{basic} so that
$Z_{n,j}\overset{a.s.}\longrightarrow Z_{(j)}$ for some random
variable $Z_{(j)}$. Define
\begin{gather*}
C_{n,j}=\sqrt{n}\,\bigl(\,\frac{1}{n}\sum_{k=1}^nX_{k,j}\,-\,Z_{n,j}\bigr)\quad\text{and}\quad
D_{n,j}=\sqrt{n}\,\bigl(Z_{n,j}-Z_{(j)}\bigr).
\end{gather*}
Next result is quite expected at this point.
\begin{cor}\label{priophg}
Suppose conditions ($**$), \eqref{hg5dr}, \eqref{pawbn} hold and fix
$1\leq j\leq d$. Then,
\begin{gather*}
\bigl(C_{n,j},\,D_{n,j}\bigr)\longrightarrow\mathcal{N}(0,U_j)\times\mathcal{N}(0,V_j)\quad\text{stably,
where}
\\U_j=V_j-Z_{(j)}(1-Z_{(j)})\quad\text{and}\quad
V_j=\frac{Z_{(j)}}{m^2}\,\bigl\{\,q_j\,(1-Z_{(j)})^2\,+\,Z_{(j)}\sum_{i\neq
j}q_i\,Z_{(i)}\,\bigr\}.
\end{gather*}
Moreover,
$E\bigl(f(D_{n,j})\mid\mathcal{G}_n\bigr)\overset{a.s.}\longrightarrow\mathcal{N}(0,V_j)(f)$
for each $f\in C_b(\mathbb{R})$, that is, $D_{n,j}$ meets condition
\eqref{pevb}.
\end{cor}
\begin{proof}
Just repeat the proof of Corollary \ref{poi} with $X_{n,j}$ in the place of $X_n$.
\end{proof}

A vectorial version of Corollary \ref{priophg} can be obtained with
slight effort. Let $\mathcal{N}_d(0,\Sigma)$ denote the
$d$-dimensional Gaussian law with mean vector 0 and covariance
matrix $\Sigma$ and
\begin{equation*}
{\bf C_n}=\bigl(C_{n,1},\ldots,C_{n,d}\bigr),\quad {\bf
D_n}=\bigl(D_{n,1},\ldots,D_{n,d}\bigr).
\end{equation*}

\begin{cor}\label{prio98poi}
Suppose conditions ($**$), \eqref{hg5dr}, \eqref{pawbn} hold. Then,
\begin{gather*}
\bigl({\bf C_n},\,{\bf D_n}\bigr)\longrightarrow\mathcal{N}_d(0,{\bf
U})\times\mathcal{N}_d(0,{\bf V})\quad\text{stably,}
\end{gather*}
where ${\bf U},\,{\bf V}$ are the $d\times d$ matrices with entries
$U_{j,j}=U_j$, $V_{j,j}=V_j$, and
\begin{gather*}
U_{i,j}=V_{i,j}+Z_{(i)}Z_{(j)},\quad
V_{i,j}=\frac{Z_{(i)}Z_{(j)}}{m^2}\,\bigl\{\sum_{h=1}^dq_hZ_{(h)}-q_i-q_j\bigr\}\quad\text{
for }i\neq j.
\end{gather*}
Moreover, $E\bigl(f({\bf
D_n})\mid\mathcal{G}_n\bigr)\overset{a.s.}\longrightarrow\mathcal{N}_d(0,{\bf
V})(f)$ for each $f\in C_b(\mathbb{R}^d)$.
\end{cor}
\begin{proof}
Given a linear functional $\phi:\mathbb{R}^d\rightarrow\mathbb{R}$,
it suffices to see that
\begin{gather*}
\phi({\bf C_n})\longrightarrow\mathcal{N}_d(0,{\bf
U})\circ\phi^{-1}\quad\text{stably, and} \\E\bigl(g\circ \phi({\bf
D_n})\mid\mathcal{G}_n\bigr)\overset{a.s.}\longrightarrow\mathcal{N}_d(0,{\bf
V})(g\circ \phi)\quad\text{for each }g\in C_b(\mathbb{R}).
\end{gather*}
To this purpose, note that
\begin{gather*}
\phi({\bf
C_n})=\sqrt{n}\,\bigl\{\,\frac{1}{n}\sum_{k=1}^n\phi(X_{k,1},\ldots,X_{k,d})\,-\,E\bigl(\phi(X_{n+1,1},\ldots,X_{n+1,d})\mid\mathcal{G}_n\bigr)\,\bigr\},
\\\phi({\bf
D_n})=\sqrt{n}\,\bigl\{\,E\bigl(\phi(X_{n+1,1},\ldots,X_{n+1,d})\mid\mathcal{G}_n\bigr)\,-\,\phi(Z_{(1)},\ldots,Z_{(d)})\,\bigr\},
\end{gather*}
and repeat again the proof of Corollary \ref{poi} with
$\phi(X_{n,1},\ldots,X_{n,d})$ in the place of $X_n$.
\end{proof}

A nice consequence of Corollary \ref{prio98poi} is that
\begin{equation*}
{\bf W_n}={\bf C_n+D_n}\longrightarrow\mathcal{N}_d(0,{\bf
U+V})\quad\text{stably}
\end{equation*}
provided conditions ($**$)-\eqref{hg5dr}-\eqref{pawbn} hold, where
${\bf W_n}=\bigl(W_{n,1},\ldots,W_{n,d}\bigr)$ and \linebreak
$W_{n,j}=\sqrt{n}\,\bigl(\,\frac{1}{n}\sum_{k=1}^nX_{k,j}\,-\,Z_{(j)}\bigr)$.

Finally, we briefly mention a possible development of the above
material. Suppose condition \eqref{pawbn} is turned into
\begin{gather*}
EA_{n,j}=EA_{n,1}\quad\text{whenever }n\geq 1\text{ and }1\leq j\leq
d_0,
\\\liminf_n\,\bigl(EA_{n,1}-EA_{n,j}\bigr)>0\quad\text{whenever }j>d_0,\notag
\\m:=\lim_nEA_{n,1}>0,\quad q_j:=\lim_nEA_{n,j}^2\quad\text{whenever }1\leq j\leq d_0,\notag
\end{gather*}
for some integer $1\leq d_0\leq d$. Roughly speaking, this means
that some colors (those labelled from $d_0+1$ to $d$) are dominated
by the others. So far, we dealt with $d_0=d$ but the case $d_0<d$ is
not unusual in applications. The main trouble is that condition
\eqref{aw38j} may fail when $d_0<d$. It is still possible to get a
CLT but one should decide how to handle dominated colors. There are
essentially two options.

One is to make assumptions on dominated colors. A classical
assumption is
\begin{equation*}
\limsup_n\frac{EA_{n,j}}{EA_{n,1}}<\frac{1}{2}\quad\text{for each
}j>d_0.
\end{equation*}
Under this condition, using some ideas from \cite{MF}, an analogous
of Corollary \ref{priophg} can be proved for $(C_{n,j},\,D_{n,j})$
with $j=1,\ldots,d_0$.

The other option is to neglect dominated colors, that is, to replace
$Z_{n,j}$ and $\frac{1}{n}\sum_{k=1}^nX_{k,j}$ by
\begin{gather*}
Z_{n,j}^*=\frac{a_j+\sum_{k=1}^nA_{k,j}X_{k,j}}{\sum_{i=1}^{d_0}a_i+\sum_{k=1}^n\sum_{i=1}^{d_0}A_{k,i}X_{k,i}}\quad\text{and}\quad
M_{n,j}^*=\frac{\sum_{k=1}^nX_{k,j}}{1+\sum_{k=1}^n\sum_{i=1}^{d_0}X_{k,i}}.
\end{gather*}
Again, an analogous of Corollary \ref{priophg} can be shown for
\begin{gather*}
C_{n,j}^*=\sqrt{n}\,\bigl(M_{n,j}^*-Z_{n,j}^*\bigr)\quad\text{and}\quad
D_{n,j}^*=\sqrt{n}\,\bigl(Z_{n,j}^*-Z_{(j)}\bigr),\,\,\,\,j=1,\ldots,d_0.
\end{gather*}

The case $d_0<d$ will be deepened in a forthcoming paper.

\vspace{1cm}

\begin{center}{\bf APPENDIX}\end{center}

\vspace{0.3cm}

In the notation of Subsection \ref{2col}, let
$S_n=b+r+\sum_{k=1}^n\bigl(B_kX_k+R_k(1-X_k)\bigr)$.

\begin{lem}\label{jrblk} Under the assumptions of Corollary
\ref{poi},
\begin{equation*}
\frac{n}{S_n}\longrightarrow\frac{1}{m}\quad\text{a.s. and in
}L_p\text{ for all }p>0.
\end{equation*}
\end{lem}
\begin{proof}
Let $Y_n=B_nX_n+R_n(1-X_n)$. By ($*$) and $EB_{n+1}=ER_{n+1}$,
\begin{gather*}
E\bigl(Y_{n+1}\mid\mathcal{G}_n\bigr)=EB_{n+1}\,E\bigl(X_{n+1}\mid\mathcal{G}_n\bigr)\,+\,ER_{n+1}\,E\bigl(1-X_{n+1}\mid\mathcal{G}_n\bigr)
\\=Z_n\,EB_{n+1}\,+\,(1-Z_n)\,EB_{n+1}=EB_{n+1}\overset{a.s.}\longrightarrow
m.
\end{gather*}
Since $m>0$, Lemma \ref{hftuimn} implies
$\frac{n}{S_n}=\frac{1}{S_n/n}\overset{a.s.}\longrightarrow\frac{1}{m}$.
To conclude the proof, it suffices to see that
$E(S_n^{-p})=\,$O$(n^{-p})$ for all $p>0$. Given $c>0$, define
\begin{equation*}
S_n^{(c)}=\sum_{k=1}^n\bigl\{X_k\bigl(B_k\wedge c-E(B_k\wedge c)\bigr)+(1-X_k)\bigl(R_k\wedge c-E(R_k\wedge c)\bigr)\bigr\}.
\end{equation*}
By a classical martingale inequality (see e.g. Lemma 1.5 of \cite{LT})
\begin{equation*}
P\bigl(\abs{S_n^{(c)}}>x\bigr)\leq 2\,\exp{\bigl(-x^2/2\,\,c^2\,n\bigr)}\quad\text{for all }x>0.
\end{equation*}
Since $EB_n=ER_n\longrightarrow m$ and both $(B_n)$, $(R_n)$ are uniformly integrable (as $\sup_n\,\bigl(EB_n^2+ER_n^2\bigr)<\infty$), there are $c>0$ and an integer $n_0$ such that
\begin{equation*}
m_n:=\sum_{k=1}^n\min\bigl\{E(B_k\wedge c),\,E(R_k\wedge c)\bigr\}>n\frac{m}{2}\quad\text{for all }n\geq n_0.
\end{equation*}
Fix one such $c>0$ and let $l=m/4>0$. For every $p>0$, one can write
\begin{gather*}
E(S_n^{-p})=p\,\int_{b+r}^\infty t^{-p-1}P(S_n<t)\,dt
\\\leq\frac{p}{(b+r)^{p+1}}\,\int_{b+r}^{b+r+n\,l}P(S_n<t)\,dt\,+\,p\,\int_{b+r+n\,l}^\infty t^{-p-1}\,dt.
\end{gather*}
Clearly, $p\,\int_{b+r+n\,l}^\infty t^{-p-1}\,dt=(b+r+n\,l)^{-p}=\,$O$(n^{-p})$. Further, for each $n\geq n_0$ and $t<b+r+n\,l$, since $m_n>n\,2\,l$ one obtains
\begin{gather*}
P(S_n<t)\leq P\bigl(S_n^{(c)}<t-b-r-m_n\bigr)\leq P\bigl(S_n^{(c)}<t-b-r-n\,2\,l\bigr)
\\\leq P\bigl(\abs{S_n^{(c)}}>b+r+n\,2\,l-t\bigr)\leq 2\,\exp{\bigl(-(b+r+n\,2\,l-t)^2/2\,\,c^2\,n\bigr)}.
\end{gather*}
Hence, $\int_{b+r}^{b+r+n\,l}P(S_n<t)\,dt\leq n\,2\,l\,\exp{\bigl(-n\,\frac{l^2}{2\,c^2}\bigr)}$ for every $n\geq n_0$, so that $E(S_n^{-p})=\,$O$(n^{-p})$.
\end{proof}

\begin{rem}\label{chissa1}
As in Subsection \ref{multi}, let
$S_n=\sum_{i=1}^da_i+\sum_{k=1}^n\sum_{i=1}^dA_{k,i}X_{k,i}$. Under
conditions ($**$)-\eqref{hg5dr}-\eqref{pawbn}, the previous proof
still applies to such $S_n$. Thus,
$\frac{n}{S_n}\longrightarrow\frac{1}{m}$ a.s. and in $L_p$ for all
$p>0$.
\end{rem}

\begin{proof} [{\bf Proof of Corollary \ref{poi}}] By Lemma \ref{plm}, it is enough to prove
$C_n\rightarrow\mathcal{N}(0,U)$ stably and $D_n$ meets condition
\eqref{pevb}. Recall from Subsection \ref{2col} that
\begin{gather*}
Z_{n+1}-Z_n=\frac{(1-Z_n)\,X_{n+1}\,B_{n+1}\,-\,Z_n\,(1-X_{n+1})\,R_{n+1}}{S_{n+1}}
\\\text{and }\,\,E\bigl\{\abs{E(Z_{n+1}\mid\mathcal{G}_n)-Z_n}^p\bigr\}=\,\text{O}(n^{-2p})\quad\text{for all }p>0.
\end{gather*}
In particular, condition \eqref{basic} holds and
$\sum_k\sqrt{k}\,E\Abs{\,E(Z_k\mid\mathcal{G}_{k-1})-Z_{k-1}}<\infty$.

\vspace{0.2cm}

{\bf "$D_n$ meets condition \eqref{pevb}".} By \eqref{momcon} and
Lemma \ref{jrblk},
\begin{gather*}
E\bigl\{\abs{Z_{k-1}-Z_k}^u\bigr\} \leq
E\bigl\{\frac{(B_k+R_k)^u}{S_{k-1}^u}\bigr\}=E\bigl\{(B_k+R_k)^u\bigr\}\,E(S_{k-1}^{-u})=\text{O}(k^{-u}).
\end{gather*}
Thus,
$E\bigl\{\sup_k\sqrt{k}\,\abs{Z_{k-1}-Z_k}\bigr\}^u\leq\sum_k\,k^{\frac{u}{2}}E\bigl\{\abs{Z_{k-1}-Z_k}^u\bigr\}<\infty$
as $u>2$. In view of Remark \ref{stabqc}, it remains only to prove
that
\begin{gather*}
n\sum_{k\geq n}(Z_{k-1}-Z_k)^2=n\sum_{k\geq
n}\bigl(\frac{(1-Z_{k-1})X_kB_k}{S_k}-\frac{Z_{k-1}(1-X_k)R_k}{S_k}\bigr)^2
\\=n\sum_{k\geq n}\frac{(1-Z_{k-1})^2X_kB_k^2}{(S_{k-1}+B_k)^2}\,+\,n\sum_{k\geq
n}\frac{Z_{k-1}^2(1-X_k)R_k^2}{(S_{k-1}+R_k)^2}
\end{gather*}
converges a.s. to $V=Z(1-Z)\,\frac{(1-Z)q+Zs}{m^2}$. It is enough to
show that
\begin{equation*}
n\sum_{k\geq
n}\frac{(1-Z_{k-1})^2X_kB_k^2}{(S_{k-1}+B_k)^2}\overset{a.s.}\longrightarrow
Z(1-Z)^2\frac{q}{m^2}\,\text{ and }\, n\sum_{k\geq
n}\frac{Z_{k-1}^2(1-X_k)R_k^2}{(S_{k-1}+R_k)^2}\overset{a.s.}\longrightarrow
Z^2(1-Z)\frac{s}{m^2}\,.
\end{equation*}
These two limit relations can be proved by exactly the same
argument, and thus we just prove the first one. Let
$U_n=B_nI_{\{B_n\leq \sqrt{n}\}}$. Since $P(B_n>\sqrt{n})\leq
n^{\frac{-u}{2}}EB_n^u$, condition \eqref{momcon} yields $P(B_n\neq
U_n,$ i.o.$)=0$. Hence, it suffices to show that
\begin{equation}\label{vhqosd}
n\sum_{k\geq
n}\frac{(1-Z_{k-1})^2X_kU_k^2}{(S_{k-1}+U_k)^2}\overset{a.s.}\longrightarrow
Z(1-Z)^2\frac{q}{m^2}.
\end{equation}
Let $Y_n=n^2\frac{(1-Z_{n-1})^2X_nU_n^2}{(S_{n-1}+U_n)^2}$. Since
$(B_n^2)$ is uniformly integrable, $EU_n^2\longrightarrow q$.
Furthermore, $\frac{S_n}{n}\overset{a.s.}\longrightarrow m$ and
$Z_n\overset{a.s.}\longrightarrow Z$. Thus,
\begin{gather*}
E\bigl(Y_{n+1}\mid\mathcal{G}_n\bigr)\leq
(1-Z_n)^2(n+1)^2E\bigl(\frac{X_{n+1}U_{n+1}^2}{S_n^2}\mid\mathcal{G}_n\bigr)
\\=Z_n(1-Z_n)^2\frac{(n+1)^2}{S_n^2}\,EU_{n+1}^2\overset{a.s.}\longrightarrow
Z(1-Z)^2\frac{q}{m^2}\quad\text{ and}
\\E\bigl(Y_{n+1}\mid\mathcal{G}_n\bigr)\geq
(1-Z_n)^2(n+1)^2E\bigl(\frac{X_{n+1}U_{n+1}^2}{(S_n+\sqrt{n+1})^2}\mid\mathcal{G}_n\bigr)
\\=Z_n(1-Z_n)^2\frac{(n+1)^2}{(S_n+\sqrt{n+1})^2}\,EU_{n+1}^2\overset{a.s.}\longrightarrow
Z(1-Z)^2\frac{q}{m^2}.
\end{gather*}
By Lemma \ref{hftuimn}, for getting relation \eqref{vhqosd}, it
suffices that $\sum_n\frac{EY_n^2}{n^2}<\infty$. Since
\begin{gather*}
\frac{EU_n^4}{n^2}\leq\frac{E\bigl\{B_n^2I_{\{B_n^2\leq\sqrt{n}\}}\bigr\}}{n^{\frac{3}{2}}}+\frac{E\bigl\{B_n^2I_{\{B_n^2>\sqrt{n}\}}\bigr\}}{n}
\leq\frac{EB_n^2}{n^{\frac{3}{2}}}+\frac{EB_n^u}{n^{1+\frac{u-2}{4}}}\,,
\end{gather*}
condition \eqref{momcon} implies $\sum_n\frac{EU_n^4}{n^2}<\infty$.
By Lemma \ref{jrblk}, $E(S_{n-1}^{-4})=\,$O$(n^{-4})$. Then,
\begin{gather*}
\sum_n\frac{EY_n^2}{n^2}\leq
\sum_nn^2E\bigl\{\frac{U_n^4}{S_{n-1}^4}\bigr\}=\sum_nn^2E(S_{n-1}^{-4})\,EU_n^4\leq
c\,\sum_n\frac{EU_n^4}{n^2}<\infty
\end{gather*}
for some constant $c$. Hence, condition \eqref{vhqosd} holds.

\vspace{0.2cm}

{\bf "$C_n\rightarrow\mathcal{N}(0,U)$ stably".} By Theorem
\ref{main}, it suffices to check conditions (a) and (b) with
$U=Z(1-Z)\,\bigl(\frac{(1-Z)q+Zs}{m^2}-1\bigr)$. As to (a), since
$E\bigl\{\abs{Z_{k-1}-Z_k}^u\bigr\}=\,$O$(k^{-u})$,
\begin{gather*}
\bigl(\,n^{-\frac{1}{2}}\,E\bigl\{\max_{1\leq k\leq
n}k\,\abs{Z_{k-1}-Z_k}\bigr\}\,\bigr)^u\leq
n^{-\frac{u}{2}}\,\sum_{k=1}^nk^uE\bigl\{\abs{Z_{k-1}-Z_k}^u\bigr\}\longrightarrow
0.
\end{gather*}
We next prove condition (b). After some algebra, one obtains
\begin{gather*}
E\bigl\{(X_n-Z_{n-1})(Z_{n-1}-Z_n)\mid\mathcal{G}_{n-1}\bigr\}
=-Z_{n-1}(1-Z_{n-1})\,E\bigl\{\frac{B_n}{S_{n-1}+B_n}\mid\mathcal{G}_{n-1}\bigr\}\,+\\+\,Z_{n-1}^2(1-Z_{n-1})\,E\bigl\{\frac{B_n}{S_{n-1}+B_n}-\frac{R_n}{S_{n-1}+R_n}\mid\mathcal{G}_{n-1}\bigr\}\quad\text{a.s.}.
\end{gather*}
Arguing as in the first part of this proof ("$D_n$ meets condition
\eqref{pevb}"),
\begin{equation*}
n\,E\bigl\{\frac{B_n}{S_{n-1}+B_n}\mid\mathcal{G}_{n-1}\bigr\}\overset{a.s.}\longrightarrow
1\quad\text{and}\quad
n\,E\bigl\{\frac{R_n}{S_{n-1}+R_n}\mid\mathcal{G}_{n-1}\bigr\}\overset{a.s.}\longrightarrow
1.
\end{equation*}
Thus,
$n\,E\bigl\{(X_n-Z_{n-1})(Z_{n-1}-Z_n)\mid\mathcal{G}_{n-1}\bigr\}\overset{a.s.}\longrightarrow
-Z(1-Z)$. Further,
\begin{equation*}
E\bigl\{\bigl(X_n-Z_{n-1})^2\mid\mathcal{G}_{n-1}\bigr\}=Z_{n-1}-Z_{n-1}^2\overset{a.s.}\longrightarrow
Z(1-Z). \end{equation*} Thus, Lemma \ref{hftuimn} implies
\begin{gather*}
\frac{1}{n}\sum_{k=1}^n(X_k-Z_{k-1})^2\,+\,\frac{2}{n}\sum_{k=1}^nk\,(X_k-Z_{k-1})\,(Z_{k-1}-Z_k)\overset{a.s.}\longrightarrow
-Z(1-Z).
\end{gather*}
Finally, write
$\frac{1}{n}\sum_{k=1}^nk^2(Z_{k-1}-Z_k)^2=\frac{1}{n}\sum_{k=1}^nk^2\bigl\{\,\frac{(1-Z_{k-1})^2X_kB_k^2}{(S_{k-1}+B_k)^2}+\frac{Z_{k-1}^2(1-X_k)R_k^2}{(S_{k-1}+R_k)^2}\,\bigr\}$.
By Lemma \ref{hftuimn} and the same truncation technique used in the
first part of this proof,
$\frac{1}{n}\sum_{k=1}^nk^2(Z_{k-1}-Z_k)^2\overset{a.s.}\longrightarrow
V$. Squaring,
\begin{gather*}
\frac{1}{n}\sum_{k=1}^n\bigr\{X_k-Z_{k-1}+k(Z_{k-1}-Z_k)\bigl\}^2\,\overset{a.s.}\longrightarrow
V-Z(1-Z)=U,
\end{gather*}
that is, condition (b) holds. This concludes the proof.
\end{proof}

\end{document}